\documentclass[12pt]{article}
\usepackage{amsmath}
\usepackage{amssymb}

\usepackage[dvips]{graphicx}
\usepackage[usenames,dvipsnames,svgnames]{xcolor}

\newtheorem{theorem}{Theorem}[section]
\newtheorem{definition}[theorem]{Definition}
\newtheorem{lemma}[theorem]{Lemma}

\newenvironment{proof}
{\noindent
{\bf Proof.}}
{\hfill $\square$\medskip

}

\newcommand{\R}{\mathbb{R}}

\newcommand{\N}{\mathbb{N}}
\renewcommand{\P}{\mathbb{P}}

\newcommand{\bv}{{\mathbf v}}

\newcommand{\bs}{{\mathbf s}}
\newcommand{\bx}{{\mathbf x}}
\newcommand{\by}{{\mathbf y}}

\newcommand{\bn}{{\mathbf n}}

\begin{document}
\title{Sobolev Regularity of Isogeometric Finite 
Element Spaces with Degenerate Geometry Map}

\author{Ulrich Reif}

\date{\today}
\maketitle

\begin{abstract}
We investigate Sobolev regularity
of bivariate functions obtained 
in Isogeometric Analysis when using geometry maps
that are degenerate in the sense that the first partial 
derivatives vanish at isolated points. In particular, we
show how the known $C^1$-conditions for D-patches have to be 
tightened 
to guarantee square integrability of second partial
derivatives, as required when computing finite element 
approximations of elliptic fourth order PDEs like the 
biharmonic equation.
\end{abstract}

\section{Introduction}
To approximate the solution of an elliptic partial differential 
equation of order $2k$ on the domain $\Omega$ in the sense of finite 
elements, one derives a corresponding energy functional 
$ J : W^{k,p}(\Omega) \to \R$, chooses a finite-dimensional
affine subspace $F \subset W^{k,p}(\Omega)$, and determines
the minimizer of $J$ among all $f \in F$. Most applications
require $p=2$, but also other exponents $p$
may be relevant when solving nonlinear problems like the 
$p$-Laplace equation $(k=1)$, see \cite{Barrett:1993}, or the 
$p$-biharmonic equation $(k=2)$, see \cite{Pryer:2019}.
In the following, we address finite element subspaces 
$F \subset W^{k,p}(\Omega)$ 
for orders $k \in \{1,2\}$, exponents $p \in [1,\infty)$,
and planar domains $\Omega\subset \R^2$.

Sobolev regularity of the shape functions spanning $F$ is usually
easy to determine if they are defined directly on a 
tessellation of $\Omega$. For instance, the ubiquitous Courant
element is $W^{1,\infty}$, while piecewise polynomial $C^1$-elements 
like Clough-Tocher, Powell-Sabin, Argyris, etc.\ are $W^{2,\infty}$.
However,
in Isogeometric Analysis \cite{Hughes:2009}, the situation is 
potentially more complicated. Given an auxiliary domain $\Gamma = 
[0,1]^2 \times \{1,\dots,N\}$, consisting of a finite number of 
copies of the unit square, piecewise polynomial functions $b_i : 
\Gamma \to \R$ are used to represent both a fixed parametrization 
$\bx = \sum_i \bx_i b_i: \Gamma \to \Omega$ of the domain, also 
called {\em geometry map}, and functions
$z = \sum_i z_i b_i : \Gamma \to \R$.
Together, they define the elements of $F$ by 
$f = z \circ \bx^{-1}$.
Now it may occur that $f$ does not inherit 
full regularity from its constituents $\bx$ and $z$
if the geometry map
is degenerate in the sense that its Jacobian determinant
$\det D\bx$ has zeros on $\Gamma$. In particular, $b_i$ being 
$C^1$ and piecewise polynomial does not imply 
$W^{2,p}$-regularity of functions $ f$ in this case.

To understand regularity of 
$f = z \circ \bx^{-1}$, one observes that its graph equals
the trace of the parametrized surface%
\footnote{ 
Here and below, we use comma-separated lists in parens to 
denote elements in $\R^d$, formally understood as column vectors.
Square brackets are used to denote vectors and matrices 
of dimension as printed.}
$\bs  := (\bx, z) = \sum_i (\bx_i, z_i) b_i : \Gamma \to \R^3$, 
i.e.,
\[
   \bigl\{(x,y,f(x,y)) : (x,y) \in \Omega\bigr\} =
   \bigl\{\bs(u,v,j) : (u,v,j) \in \Gamma\bigr\}
   .
\]
Hence, the analytic smoothness of $f$ is closely related 
to the geometric smoothness of the surface $\bs$. 
In the case $k=p=2$, as occurring for linear fourth order problems 
like the biharmonic equation, 
the setting is commonly chosen such that $f$ is at least $C^1$, 
even if, strictly speaking, this is not necessary.%
\footnote{
In the bivariate case, continuity of the gradient is by no means 
necessary for square integrability of second partials. Consider,
for instance, the continuous function 
\[
 f(x,y) := \int_0^x \ln |\ln t|\, dt
 ,
\]
defined on the triangle 
$\Delta := \{(x,y) : 0 \le y \le x \le 1/2\}$.
The first partial derivative
$f_x(x,y) = \ln |\ln x|$ is unbounded, while 
$f_{xy} = f_{yy} = 0$ and
\[
  \int_\Delta f_{xx}^2 =
  \int_{x=0}^{1/2} \int_{y=0}^x \frac{dx dy}{x^2 \ln^2(x)} =
  \int_{x=0}^{1/2} \frac{dx}{x \ln^2(x)} = 
  \frac{1}{\ln 2}
  .
\]
This peculiar example, illustrating 
the breakdown of the embedding $W^{2,p} \hookrightarrow C^1$ 
in the limiting case $p=2$, should not be understood as a
stimulus to construct shape functions of that kind.
For all practical purposes, it is sensible to assume $C^1$
when $W^{2,2}$ is needed.}
Equivalently, this means that the surface $\bs$ is $G^1$
in the sense that is possesses a continuous field of 
unit normal vectors.

Based on this insight, there 
is a new interest in methods developed formerly in Computer Aided 
Design for the construction of smooth surfaces of arbitrary topology.
In this report, we focus on the analysis of such an approach 
based on geometry maps with singularities of the form
\begin{equation} 
\label{eq:dx=0}
  \partial_1 \bx(0,0,j) = \partial_2 \bx(0,0,j) = 0
\end{equation}
for certain indices $j$.
That is, these patches of $\bx$ have vanishing first 
partial derivatives at one of their corners, which w.l.o.g.\ is
assumed to be the origin. For surface design, singular 
parametrizations of that type were first suggested in 
\cite{Peters:1991}
and later on established rigorously in \cite{Reif:1994}.

In the latter reference, a special class of singularly parametrized 
surfaces, called {\em D-patches} (the {\em D} stands for 
{\em degenerate}), is introduced that enables the construction of 
geometrically smooth surfaces. Hence, the corresponding
functions $f = z \circ \bx^{-1}$ considered here are $C^1$.
The construction is simple, the function spaces 
have a natural refinement property, and the generalization 
to the important 3d case seems to be practicable 
\cite{Peters:2020, Youngquist:2022}.
By contrast, trivariate analogues of the competing 
approaches {\em subdivision} and {\em geometric continuity}
appear to be much harder to develop, see \cite{Dietz:2023} for a 
discussion of this issue.

Therefore, the concept of singular parametrization 
according to \eqref{eq:dx=0}
has attracted the interest of researchers in Isogeometric Analysis 
\cite{Nguyen:2016, Speleers:2017, 
Casquero:2020, Youngquist:2022, Wei:2022}.
The numerical results reported there suggest good convergence
properties for the Laplace and other second order equations as well 
as for the biharmonic equation, which is of fourth order. However, 
caution is advised in the assessment of such empirical findings
lacking a solid theoretical foundation.
The problem coming along with the biharmonic equation $(k=p=2)$ is 
that continuous differentiability of $f$, as ensured by the D-patch
setting, does not imply the required square integrability of second 
partials. And in fact, we will show that the corresponding functions 
$f$ are in $W^{2,p}$ for $p \in [1,2)$, but generally not in 
$W^{2,2}$.
This marginal yet crucial lack of regularity renders the approach 
ineligible for fourth order problems unless additional, hitherto 
unknown restrictions are observed.

A first systematic study of 
Sobolev regularity of shape functions under the special conditions of 
singular parametrization was presented in \cite{Takacs:2011}.
In that paper, three different types of degeneracies
(edge collapse, collinear and vanishing partials 
at patch corners) are analyzed concerning 
$W^{1,2}$-regularity. It is found that the latter case, which is 
also the topic of this report, is not problematic in that 
respect: For geometry maps $\bx = \sum_i \bx_i b_i$ according to 
\eqref{eq:dx=0} and piecewise polynomial, but otherwise arbitrary 
$z = \sum_i z_i b_i$ the 
composition $f = z \circ \bx^{-1}$ is typically $W^{1,2}$.
A similar study concerning $W^{2,2}$-regularity is 
conducted in \cite{Juttler:2012}. However, 
the case~\eqref{eq:dx=0} of vanishing partials is not 
addressed there.

In this report, we develop a comprehensive characterization of 
$W^{1,p}$- and $W^{2,p}$-regularity of functions
$f = z \circ \bx^{-1}$ with degenerate geometry map $\bx$ according 
to \eqref{eq:dx=0} proceeding as follows:
In the next section, we fix the setting and introduce notions and 
notations. In Section~\ref{sec:aux}, we provide a series of 
auxiliary results, which are needed to establish the 
desired classification
in Section~\ref{sec:result}. The conclusion summarizes our 
findings and outlines directions for future research.

\section{D-maps and D-functions}
\label{sec:D-map}

For measurable $A \subset \R^2$, $p \in [1,\infty)$, and a 
nonnegative weight function $\mu : A \to \R_{\ge 0}$, denote by 
$L^p(A)$ and $L^p(A,\mu)$ the (standard and weighted) Lebesgue
spaces of functions $f : A \to \R^d$ with finite norm 
\[
   ||f||_p : = \left(\int_A |f|^p\right)^{1/p}
   \quad\text{and}\quad 
   ||f||_{p,\mu} : = \left(\int_A |f|^p \, \mu\right)^{1/p}
   ,
\]
respectively, where $|\cdot|$ is (for instance)
the maximum norm on $\R^d$.
With the row vectors
\[
 Df := [\partial_1 f,\, \partial_2 f]
 \quad\text{and}\quad
 D^2f := 
 [\partial^2_1 f,\, \partial_1\partial_2 f,\,
 \partial_2\partial_1f,\, \partial^2_2f]
\]
of weak first and second partial derivatives 
of a real-valued function $f : A \to \R$,
the {\em Sobolev spaces} $W^{k,p}(A)$ are defined by%
\footnote{ 
This definition is convenient for our purposes and
equivalent to the standard one based on the norm 
$\bigl(\sum_{|\alpha| \le k} \|D^\alpha f\|_p^p\bigr)^{1/p}$.
}
\[
 W^{k,p}(A) := 
 \{f : A \to \R : \|f\|_p + \|D^k f\|_{p} < \infty\}
 ,\quad 
 k \in\{1,2\}
 .
\]

Let $\bx : \Gamma \to \Omega$ be a parametrization of $\Omega \subset 
\R^2$ over the auxiliary domain 
$\Gamma = [0,1]^2 \times \{1,\dots,N\}$, chosen
such that the partial images $\Omega_j := \bx([0,1]^2,j)$
form an essentially disjoint cover 
$\bigcup_{j=1}^N \Omega_j = \Omega$ of the domain. Then,
if the function $f = z \circ \bx^{-1} : \Omega \to \R$
is $C^{k-1}$, it holds
\[
 \|D^k f\|_p^p = \int_\Omega |D^k f|^p
 =
 \sum_{j=1}^N \int_{\Omega_j} |D^k f|^p
 ,
\]
implying that $f \in W^{k,p}(\Omega)$ if and only if 
$f_j \in W^{k,p}(\Omega_j)$ for each 
segment $f_j := z \circ \bx(\cdot,\cdot,j)^{-1}$.
For this reason, we confine the following considerations to 
the analysis of geometry maps with only one segment. That is, we 
assume $\Gamma = [0,1]^2 \times \{1\}$ and omit notation of the 
insignificant index $j=1$.

Given subsets $A \subset \R^2$ and $B \subset \R^d$, we denote by
$\P(A,B)$ the space of bivariate $\R^d$-valued polynomials
\[
 \pi : A \to B
 ,\quad 
 \pi(u,v) = \sum_{j,k \ge 0} \pi_{jk} u^j v^k
 ,
\]
where only a finite number of monomial coefficients 
$\pi_{jk} \in \R^d$ is different from $0$. By convention, 
coefficients are always indicated by the same letter as the function 
itself tagged with a double subscript.

To parametrize subsets of $\R^2$ over the unit 
square, we will consider a special class of degenerate 
polynomial maps, characterized as follows:
\begin{definition}
Let $\Sigma := [0,1]^2$, 
$\Sigma^* := \Sigma \setminus \{(0,0)\}$, 
and $\Omega \subset \R^2$.
A bijective polynomial function $\bx \in \P(\Sigma,\Omega)$
is called a {\em D-map} with image $\Omega$ if
\begin{itemize}
\item[{\it i)}]
$\bx$ is degenerate at the origin according to
$\bx_{10} = \bx_{01} = \bx_{11} = 0$;
\item[{\it ii)}]
the Jacobian determinant $\det D\bx$ does not vanish 
on $\Sigma^*$;
\item[{\it iii)}]
the $(2\times 2)$-matrices $[\bx_{20},\bx_{21}],\, 
[\bx_{02},\bx_{12 }]$, and $[\bx_{20},\bx_{02}]$
have full rank.
\end{itemize}
A D-map $\by \in \P(\Sigma,\Xi)$ with image $\Xi$ 
is in {\em standard form} if
$\by_{00} = (0,0)$ is the origin and
$\by_{20} = (1,0)$,
$\by_{02} = (0,1)$ are the unit vectors in $\R^2$.
\end{definition}
%
Let us briefly motivate this definition:
The first property facilitates the 
parametrically smooth contact of $n \neq 4$
maps $\bx(\cdot,\cdot,1),\dots,\bx(\cdot,\cdot,n)$ sharing the vertex
$\bx(0,0,1) = \cdots = \bx(0,0,n)$,
while the second one admits to focus the analysis on the 
isolated singularity at the origin. 
The third property excludes degeneracies of higher order.
D-maps in standard form, for which we reserve the letters
$\by$ and $\Xi := \by(\Sigma)$
throughout, are useful to simplify the following analysis.

Isogeometric analysis employs composition
to construct approximations on the image 
of a given parametrization. Specifically, if the parametrization is 
a D-map, we consider functions of the following form:
\begin{definition}
\label{def:z}
The composition
\[
  f := z \circ \bx^{-1} : \Omega \to \R
\]
is called a {\em D-function} on $\Omega \subset \R^2$ if 
$\bx \in\P(\Sigma,\Omega)$ is a D-map with image $\Omega$ and 
$z \in \P(\Sigma,\R)$ is a real-valued polynomial.
A D-function $h := w \circ \by^{-1}$ is in {\em standard form} if 
$\by$ is in standard form and if 
the coefficients of $w \in \P(\Sigma,\R)$ satisfy 
$w_{20} = w_{02} = 0$.
\end{definition}
One may think of $\bx$ as being given and of $z$ as
to be determined suitably to serve the desired approximation 
purposes.
Special functions $\by$ and $w$ corresponding to 
$\bx$ and $z$ allow to reduce the investigation of 
D-functions to the much simpler standard form.
They are derived as follows:

First, with $T := [\bx_{20}, \bx_{02}]^{-1}$ denoting the inverse
of the $(2 \times 2)$-matrix built from the two non-zero quadratic 
coefficients of a given D-map $\bx$, we define
\[
  \by = \Phi(\bx) := T (\bx - \bx_{00})
  \in \P(\Sigma,\Xi)
  ,
\]
where $\Xi := \by(\Sigma)$. 
We call $\by$ the {\em D-map in standard form related to $\bx$}.
To show that $\by$ has all properties of a D-map
in standard form, we note that its coefficients are given by 
\[
 \by_{jk} = 
 \begin{cases}
  0 & \text{if } j=k=0, \\
  T\bx_{jk} & \text{else}.
 \end{cases}
\]
In particular, $\by_{00} = \by_{10} = \by_{01} = \by_{11} = 0$
and $\by_{20} = (1,0), \by_{02} = (0,1)$. 
Further, $\det D\by = \det T \cdot \det D\bx \neq 0$ on $\Sigma^*$ 
and the matrices 
$[\by_{20},\by_{21}] = T \cdot [\bx_{20}, \bx_{21}],\,
 [\by_{02},\by_{12}] = T \cdot [\bx_{02}, \bx_{12}],\,
 [\by_{20},\by_{02}] = T \cdot [\bx_{20}, \bx_{02}]$
have full rank.
To simplify notation, we write
\begin{equation}
\label{eq:abcd}
 \by_{21} = \begin{bmatrix} \alpha \\ \beta \end{bmatrix}
 ,\quad 
 \by_{12} = \begin{bmatrix} \gamma \\ \delta \end{bmatrix}
\end{equation}
throughout and call $\alpha,\dots,\delta$ the {\em parameters} 
of $\bx$. These are the unique coefficients satisfying
\begin{equation}
\label{eq:xcond}
 \bx_{21} = \alpha \bx_{20} + \beta \bx_{02}
 ,\quad 
 \bx_{12} = \gamma \bx_{20} + \delta \bx_{02}
 .
\end{equation}
Second, given $z \in \P(\Sigma,\R)$ and the 
D-map $\by$ in standard form, we define the polynomial 
$w \in \P(\Sigma,\R)$ by 
\[
 w = \Psi(\by,z) := z - [z_{20},\,z_{02}]\cdot \by
 .
\]
Its coefficients are 
\begin{equation}
\label{eq:w}
  w_{jk} = z_{jk} - [z_{20},\,z_{02}]\cdot \by_{jk}
  ,
\end{equation}
and in particular $w_{20}=w_{02} = 0$. With $\by = \Phi(\bx)$ 
and $w = \Psi(\by,z)$ as above, we call 
\[
  h := w \circ \by^{-1}
\]
the {\em  D-function in standard form related to 
$f = z \circ \bx^{-1}$}.

To complete the section, we establish
a relation between D-functions and D-patches 
as introduced in \cite[Definition~1.1]{Reif:1994}:
\begin{theorem}
\label{thm:C1}
Let $\bs = (\bx, z) \in \P(\Sigma,\R^3)$ be the 
parametrized surface corresponding to the D-function
$f = z \circ \bx^{-1}$. Then $\bs$ 
is a generic D-patch if
\[
 z_{10}=z_{01}=z_{11} = 0
 ,\quad 
 z_{21} = \alpha z_{20}+\beta z_{02}
 ,\quad 
 z_{12} = \gamma z_{20} + \delta z_{02}
 .
\]
In this case, $f$ is continuously differentiable.
\end{theorem}
\begin{proof}
The coefficients of $\bs$ are
$\bs_{jk} = (\bx_{jk}, z_{jk})$.
First, $\bs$ is degenerate in the sense of 
\cite{Reif:1994} because 
$\bs_{10} = \bs_{01} = \bs_{11} = 0$.
Second, the conditions of the theorem together 
with \eqref{eq:xcond} yield
\[
 \bs_{21} = \alpha \bs_{20} + \beta \bs_{02}
 ,\quad 
 \bs_{12} = \gamma \bs_{20} + \delta \bs_{02}
 .
\]
Anticipating positivity of the coefficients $\beta, \gamma$
as established in Lemma~\ref{lem:bg>0},
these are exactly the coplanarity constraints characterizing
a D-patch. Third, $\bs$ is generic since 
linear independence of $\bx_{20}$ and $\bx_{02}$ 
implies linear independence of 
$\bs_{20} = (\bx_{20}, z_{20})$ and 
$\bs_{02} = (\bx_{02}, z_{02})$.

On $\Sigma^*$, the normal of $\bs$ is 
$\bn := (\partial_1 \bs \times \partial_2 \bs)/
\|\partial_1 \bs \times \partial_2 \bs\|$.
According to \cite[Theorem 1.2]{Reif:1994},
its limit at the origin exists and is given by 
\[
 \bn_0 = (n_0^1,n_0^2,n_0^3) := \lim_{(u,v) \to (0,0)} \bn(u,v)
 =
 \frac{\bs_{20} \times \bs_{02}}{\|\bs_{20} \times \bs_{02} \|}
 ,
\]
where the third component 
$n_0^3 = \det [\bx_{20},\bx_{02}]/\|\bs_{20} \times \bs_{02} \|$ 
does not vanish. Since $\det D\bx(u,v) \neq 0$ for 
$(u,v) \neq (0,0)$, the inverse function theorem implies 
continuous differentiability of $\bx^{-1}$ and hence 
also of $f = z \circ \bx^{-1}$ 
on $\Omega^* := \Omega \setminus \{(x_0,y_0)\}$, where 
$(x_0,y_0) := \bx(0,0)$.
Since the graph of $f$ equals the trace 
of $\bs$, the normal satisfies 
\[
 \bn \circ \bx^{-1} = \frac{1}{\sqrt{1+\|\nabla f\|^2}}
 \begin{bmatrix}
  -\nabla f \cr 1
 \end{bmatrix}
\]
on $\Omega^*$.
Passing to the limit, we obtain 
\[
 \lim_{(x,y) \to (x_0,y_0)} \frac{1}{\sqrt{1+\|\nabla f(x,y)\|^2}}
 = n_0^3 \neq 0
\]
and further 
\[
   \lim_{(x,y) \to (x_0,y_0)} \nabla f(x,y) 
   =
   \frac{-(n_0^1,n_0^2)}{n_0^3}
   .
\]
Hence, $f$ is continuously differentiable on all of 
$\Omega$ with $\nabla f(x_0,y_0) = -(n_0^1, n_0^2)/n_0^3$.
\end{proof}

\section{Auxiliary results}
\label{sec:aux}

To prepare the analysis of Sobolev regularity of D-functions, 
we need to establish a series of technical results.

\begin{lemma}
\label{lem:fgh}
Given a D-map $\bx \in \P(\Sigma,\Omega)$ and a polynomial
$z \in \P(\Sigma,\R)$, let 
$\by := \Phi(\bx)$ and $w := \Psi(\by,z)$.
Then the D-functions 
\[
  f := z \circ \bx^{-1}
  ,\quad 
  g := z \circ \by^{-1}
  ,\quad 
  h := w \circ \by^{-1}
\]  
have the same Sobolev regularity, i.e.,
$f \in W^{k,p}(\Omega) \Leftrightarrow 
g \in W^{k,p}(\Xi) \Leftrightarrow 
h \in W^{k,p}(\Xi)$.
\end{lemma}
\begin{proof}
The affine map $\phi : \bv \mapsto T(\bx - \bx_{00})$, 
relating $\bx$ and $\by$ by $\by = \phi \circ \bx$,
provides a smooth diffeomorphism between $\Omega$ and $\Xi$.
Hence, $f = z \circ \bx^{-1} = z \circ \by^{-1} \circ \phi
= g \circ \phi$ have equal Sobolev regularity.
Further, defining the linear function
$\lambda : (\xi,\eta) \mapsto z_{20}\xi + z_{02}\eta$, we have
$w = z - \lambda \circ \by$ and 
$h  = w\circ \by^{-1} = g - \lambda$, implying 
equal Sobolev regularity of $g$ and $h$.
\end{proof}
In the following, it is convenient to use 
the expansion
\begin{equation}
\label{eq:xcan}
  \by (u,v) =
  \begin{bmatrix}
   u^2 + \alpha u^2v + \gamma u v^2 \\
   v^2 + \beta  u^2v + \delta u v^2
  \end{bmatrix}
  + O\big(u^3 + v^3 + (u+v)^4\big)
\end{equation}
for the D-map $\by$ in standard form with 
parameters $\alpha,\dots,\delta$ according to \eqref{eq:abcd}.
\begin{lemma}
\label{lem:iG}
The Jacobian determinant $\mu := \det D\by \in \P(\Sigma,\R)$ of 
a D-map $\by$ in standard form 
according to \eqref{eq:xcan} is given by
\begin{equation}
\label{eq:g}
  \mu(u,v) = 4uv + 2\beta u^3 + 2 \gamma v^3 + \nu(u,v)
  ,
\end{equation}
where the remainder $\nu$ is of order
\[
 \nu(u,v) = O\big(u^2v + uv^2 + (u+v)^4 \big)
 .
\]
On $\Sigma^*$, $\mu$ is positive and $D\by$ is invertible
with $D\by^{-1} = \mu^{-1}\Gamma$, where
$\Gamma \in \P(\Sigma,\R^{2\times 2})$ is the adjugate matrix of 
$D\by$, given by
\[
 \Gamma(u,v) =  
 \begin{bmatrix}
 2v+\beta u^2 +O(v^2+uv)&  O\big( (u+v)^2 \big)\\
 O\big( (u+v)^2 \big) & 2u+\gamma v^2+O(u^2+uv)
 \end{bmatrix}
 + 
 O\big( (u+v)^3 \big)
 .
\]
\end{lemma}
\begin{proof}
The formulas for $\mu$ and $\Gamma$
are easily verified by inspection. Concerning $\mu$,
we conclude from Property~{\it ii)}
that it must not vanish on $\Sigma^*$.
The leading term $4uv$ of $\mu$ is positive near
the origin, and hence, $\mu$ is positive on 
$\Sigma^*$.
\end{proof}

\begin{lemma}
\label{lem:bg>0}
The parameters $\beta$ and $\gamma$ of a D-map $\bx$ are positive.
\end{lemma}
\begin{proof}
The case $\beta=0$ is excluded by Property~{\it iii)}.
Assuming $\beta < 0$ and choosing $t>0$ sufficiently small,
it holds $(u,v) := (-\beta t^2, 2t) \in \Sigma^*$ 
and $\mu(u,v) = t^3(8\beta + O(t)) < 0$. This contradicts 
Lemma~\ref{lem:iG}.
The coefficient $\gamma$ can be treated in the same way.
\end{proof}
\begin{lemma}
\label{lem:g}
Denote the lower subtriangle of the domain $\Sigma$
by $\Delta := \{(u,v) \in \Sigma : v \le u\}$.
There exist constants $c_1,c_2>0$ such that 
the function $\mu$ given by \eqref{eq:g} is bounded by
\[
  c_1(uv + u^3) \le \mu(u,v) \le c_2(uv + u^3)
  ,\quad 
  (u,v) \in \Delta
  .
\]
\end{lemma}
\begin{proof}
On $\Delta$, the estimate $(u+v)^4 \le (2u)^3 (u+v)$
implies 
$u^2v + uv^2 + (u+v)^4 \le 8(uv+u^3)(u+v)$ and hence
$\nu(u,v) = (uv+u^3)\, O(u+v)$ for the remainder in 
\eqref{eq:g}.
Let $\varepsilon>0$ be chosen sufficiently small such that 
the term $O(u+v)$ is bounded 
by $\min\{2,\beta\}$ in modulus for 
$(u,v) \in \Delta_\varepsilon := \{(u,v) \in \Delta : 
u+v < \varepsilon\}$.
For $(u,v) \in \Delta_\varepsilon$,
\[
 \mu(u,v) \ge
 4uv + 2 \beta u^3 - |\nu(u,v)| \ge
 \min\{2,\beta\}\, (uv+u^3)
\]
and
\[
 \mu(u,v)  \le 
 4uv + 2\beta u^3 + 2 \gamma u^3 + |\nu(u,v)|
 \le (4 + 3\beta + 2\gamma)\, (uv+u^3)
 .
\]
Furthermore, the function $\rho(u,v) := \mu(u,v)/(uv+u^3)$ is 
continuous 
on the compact domain $\Delta \setminus \Delta_\varepsilon$
and, by Lemma~\ref{lem:iG}, positive. Hence, there exist
constants $c_1', c_2' > 0$ such that 
$c_1' \le \rho(u,v) \le c_2'$ for $(u,v) \in 
\Delta \setminus \Delta_\varepsilon$.
Combining the estimates and setting 
$c_1 := \min\{2,\beta,c_1'\},
 c_2 := 4+3\beta+2\gamma+c_2'$ proves the lemma.
\end{proof}

\begin{lemma}
\label{lem:int}
For exponents $m,n,\ell \in \N_0$, let 
\[
 M^{m,n}_\ell(u,v) := \frac{u^m v^n}{\mu(u,v)^\ell}
 ,\quad 
 (u,v) \in \Sigma^*
 .
\]
Then $M^{m,n}_\ell \in L^p(\Sigma,\mu)$ if and only if 
\[
  w(m,n) := m + n + \min\{m,n\} > 3\ell-6/p
  .
\]
\end{lemma}
\begin{proof}
By the previous lemma, the integrals
\[
 I := \int_{\Delta} |M^{m,n}_\ell(u,v)|^p \mu(u,v)\, dudv
 \quad\text{and}\quad
 J := 
 \int_{\Delta} \frac{u^{pm} v^{pn}}{(uv+u^3)^{p\ell-1}}\, dudv
\]
converge or diverge simultaneously.
Substituting 
$(u,v) = \varphi(s,t) := \big(st, s^2t(1-t)\big)$, 
we obtain 
$|\det D\varphi(u,v)| = s^2 t$ and 
\[
 J = 
 \int_{\varphi^{-1}(\Delta)} f_1(s) f_2(t)\, ds dt
 ,
\]
where 
\[
 f_1(s) := s^{p(m+2n-3\ell)+5}
 \quad\text{and}\quad 
 f_2(t) := t^{p(m+n-2\ell)+3} \,(1-t)^{pn}
 .
\]
Since
$\varphi\big([0,1]\times[0,1]\big) 
\subset \Delta \subset
\varphi\big([0,2] \times [0,1]\big)$, 
see Figure~\ref{fig:phi} instead of a tedious argument,
we obtain
\[
 \int_0^1 f_1(s)\, ds \int_0^1 f_2(t)\, dt 
 \le J \le 
 \int_0^2 f_1(s)\, ds \int_0^1 f_2(t)\, dt
 .
\]
That is, $J$ converges if and only if 
$m+2n>3\ell-6/p$ and $m+n>2\ell-4/p$.
\begin{figure}
\centering{
\includegraphics[height=5cm]{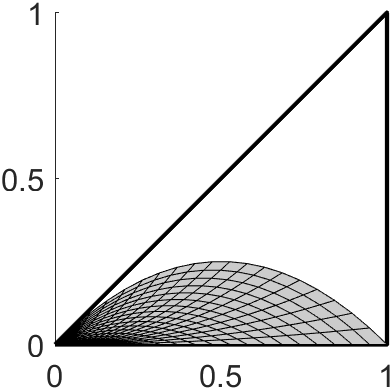}
\qquad
\includegraphics[height=5cm]{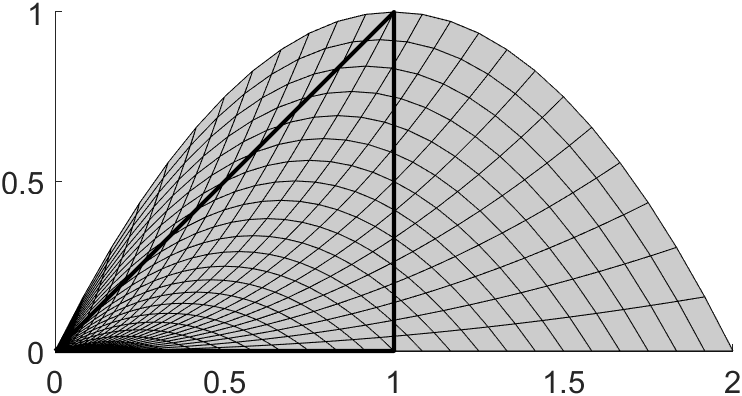}}
\caption{Image of map $\varphi$ for domains
$[0,1]\times [0,1]$ {\it (left)} and 
$[0,2]\times [0,1]$ {\it (right)}.
}
\label{fig:phi}
\end{figure}

Analogously, on the upper subtriangle
$\Delta' := \{(u,v) \in \Sigma^* : v \ge u\}$, we obtain 
integrals $I'$ and $J'$, where 
$J'$ and therewith $I'$ converges 
if and only if $2m+n>3\ell-6/p$ and $m+n>2\ell-4/p$.

Summarizing, $\|M^{m,n}_\ell\|_{p,\mu}^p = I + I'$
converges to a finite value if and only if 
\[
 m+2n>3\ell-6/p
 \quad \text{and} \quad 
 2m+n>3\ell-6/p
 \quad \text{and} \quad
 m+n > 2\ell-4/p
 .
\]
The first two inequalities are equivalent to 
$w(m,n) > 3\ell-6/p$. When added together,
they yield the third one so that this condition is 
redundant.
\end{proof}
\begin{lemma}
\label{lem:Lambda}
Let $r,\ell \in \N_0 $ and the function $\mu$ according to 
\eqref{eq:g} be given. For a polynomial
$P \in \P(\Sigma,\R^d)$,
define the {\em leading part} of order $r$ by 
\[
  \Lambda_r P := \sum_{w(j,k) \le r} P_{jk} u^j v^k.
\]
Further, let
\[
 p \in \Pi(r,\ell) := \big[6/(3\ell-r), 6/(3\ell - 1-r) \big)
 .
\]
Then it holds:
\begin{itemize}
\item[a)]
If $\Lambda_r P = 0$, then $\mu^{-\ell} P \in L^p(\Sigma,\mu)$.
\item [b)]
If there exists a component $P' \in \P(\Sigma,\R)$ of the 
vector-valued function $P$ which is a multiple of a monomial, i.e., 
$\Lambda_r P' = P'_{mn} u^m v^n$ for 
some scalar coefficient $P'_{mn}$, then
$\mu^{-\ell} P \in L^p(\Sigma,\mu)$ if and only if $\P'_{mn} = 0$.
\end{itemize}
\end{lemma}
\begin{proof}
Splitting $P$ into the leading part $\Lambda_r P$ and 
the remainder $P-\Lambda_r P$, we obtain the 
estimates 
\begin{align}
  \label{eq:le}
  \|\mu^{-\ell} P \|_{p,\mu}
  &\le 
  \|\mu^{-\ell} \Lambda_r P  \|_{p,\mu} +
  \|\mu^{-\ell} (P-\Lambda_r P)\|_{p,\mu} \\
  \label{eq:ge}
  \|\mu^{-\ell} P \|_{p,\mu}
  &\ge 
  \|\mu^{-\ell} \Lambda_r P  \|_{p,\mu} -
  \|\mu^{-\ell} (P-\Lambda_r P) \|_{p,\mu}
 .
\end{align}
The condition $p\in \Pi(r,\ell)$ is equivalent to
$r+1 > 3\ell - 6/p \ge r$.
Hence, by Lemma~\ref{lem:int},
\begin{equation}
\label{eq:T}
  \|\mu^{-\ell} (P-\Lambda_r P) \|_{p,\mu} 
  \le 
  \sum_{w(j,k) \ge r+1} 
  \|  P_{jk} M^{j,k}_\ell \|_{p,\mu}
  \le
  \sum_{w(j,k) \ge r+1} |P_{jk}|\cdot
  \| M^{j,k}_\ell \|_{p,\mu}
  < \infty
  .
\end{equation}
Concerning a), if $\Lambda_r P = 0$, then 
$\mu^{-\ell} P \in L^p(\Sigma,\mu)$
follows immediately from \eqref{eq:le} and \eqref{eq:T}.
Concerning b), assume that $\mu^{-\ell} P \in L^p(\Sigma,\mu)$.
Then, by \eqref{eq:ge} and \eqref{eq:T}, also 
$\mu^{-\ell} \Lambda_r P \in L^p(\Sigma,\mu)$.
This implies that each component of $\mu^{-\ell}\Lambda_r P$ lies 
in $L^p(\Sigma,\mu)$. In particular,
$\mu^{-\ell}\Lambda_r P' = P'_{mn} M^{m,n}_\ell
\in L^p(\Sigma,\mu)$. 
Unless $P'_{mn} = 0$, this contradicts Lemma~\ref{lem:int} 
since $w(m,n) \le r \le 3\ell-6/p$.
\end{proof}

\section{Sobolev regularity}
\label{sec:result}
In the context of finite element simulations,
Sobolev regularity of D-functions $f = z \circ \bx^{-1}$ is 
crucial. For most applications,
$f \in W^{1,2}(\Omega)$ or $f \in W^{2,2}(\Omega)$ is needed, but 
nonlinear equations like the $p$-Laplace or the 
$p$-biharmonic equation have different requirements.
Therefore, it is our goal to provide a comprehensive list of 
criteria for the membership of $f$ in $W^{k,p}(\Omega)$
depending on $k \in \{1,2\}$ and $p \in [1,\infty)$.
For $k=1$, the vanishing of certain coefficients $z_{jk}$
suffices to characterize three different cases.
By contrast, for $k=2$, the conditions for some of the six 
cases become rather involved when formulated in terms of 
the coefficients of $z$ and $\bx$. Here, it is much more
convenient to specify them for the coefficients of the functions
$w$ and $\by$ defining the related standard form 
$h = w \circ \by^{-1}$.

\begin{theorem}
\label{thm:H1}
Let $f = z \circ \bx^{-1}$ be a D-function. Then 
$f \in W^{1,p}(\Omega)$ for all 
\begin{itemize}
\item[a)] 
$p \in [1,3)$;
\item[b)]
$p \in [3,6)$ iff $z_{10} = z_{01} = 0$;
\item[c)] 
$p \in [6,\infty)$ iff $z_{10} = z_{01} = z_{11} = 0$.
\end{itemize}
\end{theorem}
\begin{proof}
Replacing $\bx \in \P(\Sigma,\Omega)$ by the related 
D-map $\by = \Phi(\bx) \in \P(\Sigma,\Xi)$ in standard form greatly 
simplifies the derivation of the result.
Concretely, by Lemma~\ref{lem:fgh}, 
$f \in W^{1,p}(\Omega)$ if and only if
$g  = z \circ \by^{-1} \in W^{1,p}(\Xi)$.
Since $g$ is a continuous function on the compact 
domain $\Xi$, finiteness of $\|g\|_p$ can be taken 
for granted so that it remains to consider finiteness
of $\|Dg\|_p$.
Differentiation of $z = g \circ \by$
yields $Dz = Dg \circ \by \cdot D\by$. 
Hence, using the parametrization $\by : \Sigma \to \Xi$
of the domain of $g$ and recalling 
$D\by^{-1} = \mu^{-1} \Gamma$ from Lemma~\ref{lem:iG},
\[
 \|Dg\|_p^p = 
 \int_{\Xi} |Dg|^p = 
 \int_{\Xi} |(Dz\, D\by^{-1}) \circ \by^{-1}|^p = 
 \int_{\Sigma} \big(\mu^{-1} |Dz \,\Gamma|\big)^p \, \mu = 
 \|\mu^{-1} Dz\, \Gamma \|_{p,\mu}^p
 .
\]
Setting $P:= Dz\,\Gamma \in \P(\Sigma,\R^2)$, finiteness of 
$\|Dg\|_p$ is equivalent to $\mu^{-1} P \in 
L^p(\Sigma,\mu)$, which can now be decided using 
Lemma~\ref{lem:Lambda} for $\ell=1$ and different values of $r$.
We consider the three cases in turn:
\begin{itemize}
\item[a)]
Let $r=0$ and $p \in \Pi(0,1) = [2,3)$.
The leading part $\Lambda_0 P = [0, 0]$ vanishes and finiteness
of $\|Dg\|_p$ follows.
This implies also finiteness for smaller exponents $p \in [1,2)$.
\item[b)] 
Let $r=1$ and $p \in \Pi(1,1) = [3,6)$.
The leading part is
$\Lambda_1 P = [2z_{10}v, 2z_{01}u]$.
Hence, $\|Dg\|_p$ is finite if $z_{10} = z_{01}=0$. 
Both components of $\Lambda_1 P$ 
are multiples of monomials, implying that $\|Dg\|_3$ is not finite
if either $z_{10} \neq 0$ or $z_{01} \neq 0$.
\item[c)] 
Let $r=2$ and $p \in \Pi(2,1) = [6,\infty)$.
In view of the preceding case, $\|Dg\|_p$ can be finite 
only if $z_{10} = z_{01}=0$. 
Assuming that, the leading part is 
$\Lambda_3 P = [2z_{11}v^2, 2z_{11}u^2]$.
Hence, $\|Dg\|_p$ is finite if $z_{11}=0$.
Both components of $\Lambda_3 P$ 
are multiples of monomials, implying that $\|Dg\|_6$ is not finite
if $z_{11} \neq 0$.
\end{itemize}
\end{proof}
The theorem states that
square integrability of the gradient of a D-function 
comes for free in accordance with the findings in 
\cite{Takacs:2011}. Further, mimicking the singularity of the 
D-map $\bx$ by the polynomial $z$ increases the range of
admissible exponents $p$. We remark that the strongest condition 
$z_{10}=z_{01}=z_{11}=0$ even yields boundedness of the 
gradient, $f \in W^{1,\infty}(\Omega)$, as can be shown using the 
substitution
$(u,v) = \varphi(s,t)$ introduced in the proof of 
Lemma~\ref{lem:int}. The following result concerning 
second partials is richer in content:

\begin{theorem}
\label{thm:H2}
Let $f = z \circ \bx^{-1}$ be a D-function
and $h = w \circ \by^{-1}$ be the related standard form. 
Then $f \in W^{2,p}(\Omega)$ for all 
\begin{itemize}
\item[a)]
$p \in [1,6/5)$    iff 
$w_{10}=w_{01} = 0$;
\item[b)]
$p \in [6/5,3/2)$  iff $w$ as above and
$w_{11} = 0$;
\item[c)]
$p \in [3/2,2)$    iff $w$ as above and
$w_{21} = w_{12} = 0$;
\item[d)]
$p \in [2,3)$      iff $w$ as above and
$w_{31} = 3\alpha w_{30}/2$,
$w_{13} = 3\delta w_{03}/2$;
\item[e)] $p \in [3,6)$      iff $w$ as above and
$w_{30} = w_{31} = w_{03} = w_{13} = 0$ and 
\[
  w_{41} = 2\alpha w_{40} + \beta w_{22}
  ,\quad 
  w_{14} = 2\delta w_{04} + \gamma w_{22}
  .
\]
\item[f)] $p \in [6,\infty)$ iff $w$ as above and
\begin{align*}
  w_{51} &=
  [w_{22}/2,\, w_{40}] \cdot
  (2\by_{31} - 3\alpha\by_{30} - 2\beta \by_{12}) + 
  \beta w_{32} + 5 \alpha w_{50}/2 \\
  w_{15} &=
  [w_{22}/2,\, w_{04}] \cdot
  (2\by_{13} - 3\delta\by_{03} - \,2\gamma \by_{21}) \,+ 
  \gamma w_{23} + 5 \delta w_{05}/2 
  .
\end{align*}
\end{itemize}
\end{theorem}
\begin{proof}
As in the preceding proof, $\|h\|_p$ is always finite so that
it suffices to analyze $\|D^2h\|_p$. 
Denoting by $\otimes$ the Kronecker product,
twofold differentiation of $w = h \circ \by$ yields 
\begin{align*}
  Dw   &= Dh \circ \by \cdot D\by \\
  D^2w &= D^2h \circ \by \cdot (D\by \otimes D\by)
  + Dh \circ \by \cdot D^2 \by
  .
\end{align*}
Solving for $D^2 h$, we obtain 
\[
  D^2 h \circ \by = 
  (D^2w - Dw \cdot D\by^{-1} \cdot D^2 \by)\cdot
  (D\by^{-1} \otimes D\by^{-1}) = \mu^{-3} P 
  ,
\]
where 
\[
  P := 
  (\mu D^2w - Dw \cdot \Gamma \cdot D^2 \by) \cdot 
  (\Gamma\otimes \Gamma)
  \in \P(\Sigma,\R^4)
  .
\]
Hence, using the parametrization $\by : \Sigma \to \Xi$
of the domain of $h$,
\[
 \|D^2h\|_p^p = 
 \int_{\Xi} |D^2 h|^p = 
 \int_{\Sigma} \big(\mu^{-3} |P|\big)^p \, \mu = 
 \big\|\mu^{-3} P\big\|_{p,\mu}^p
 .
\]
That is, finiteness of 
$\|D^2h\|_p$ is equivalent to $\mu^{-3} P \in 
L^p(\Sigma,\mu)$, which can now be decided in almost 
all cases using 
Lemma~\ref{lem:Lambda} for different values of $r$.

The six cases specified in the theorem 
correspond to choosing $r = 3,\dots,8$.
For each $r$, the argument proceeds as follows:
First, the given range of $p$ is just $\Pi(r,3)$.
Second, one computes the leading part 
$\Lambda_r P$ of $P$
under the assumption that the conditions for the preceding 
case are satisfied since, otherwise, finiteness 
of $\|D^2h\|_p$ is not possible. 
Third, one verifies that the conditions specified in the theorem 
are equivalent to $\Lambda_r P=0$ and thus guarantee 
finiteness of $\|D^2h\|_p$. 
Fourth, one shows that $\|D^2h\|_p$ is not finite 
for the minimal value $p = 6/(9-r)$ of the given range
if one of the conditions is not satisfied. 
For $r \neq 7$, this follows immediately from the fact that 
the components of $\Lambda_r P$ are multiples of monomials. 
Only case e), corresponding to $r=7$,
requires a special treatment.

Now for the details: Using a computer algebra system to perform the 
tedious computations of 
$[a_r,b_r,c_r,d_r] := \Lambda_r P/16$, we find
$b_r = c_r = 0$ throughout. Furthermore,
the other two components are given by 
\[
\begin{array}{ll}
a_3(u,v) = -w_{10} \, v^3 &
d_3(u,v) = -w_{01} \, u^3 \\
a_4(u,v) = -w_{11} \, v^4 &
d_4(u,v) = -w_{11} \, u^4 \\
a_5(u,v) = -w_{12} \, v^5 &
d_5(u,v) = -w_{21} \, u^5 \\
a_6(u,v) = (3\delta w_{03}/2-w_{13})  v^6 \quad &
d_6(u,v) = (3\alpha w_{30}/2- w_{31})  u^6 
\end{array}
\]
and, with certain constants $a_7^*, d_7^*$ depending only on the 
coefficients of $\by$,
\begin{align*} 
  a_7(u,v) &=  w_{30} u^2 v^3 + 3 \gamma w_{30} u v^5  + 
  (a_7^* w_{03} + \gamma w_{22} + 2 \delta w_{04} - w_{14})v^{7}\\
  d_7(u,v) &=  w_{03} u^3 v^2 + 3\beta w_{03} u^5 v + 
  (d_7^* w_{30} + \beta  w_{22} + 2 \alpha w_{40} - w_{41})u^{7}\\
  a_8(u,v) &= \big(
  [w_{22}/2,\, w_{04}] \cdot
  (2\by_{13} - 3\delta\by_{03} - \,2\gamma \by_{21}) \,+ 
  \gamma w_{23} + 5\delta w_{05}/2 - w_{15}
  \big) v^8 \\
  d_8(u,v) &= \big(
  [w_{22}/2,\, w_{40}] \cdot
  (2\by_{31} - 3\alpha\by_{30} - 2\beta \by_{12}) + 
  \beta w_{32} + 5 \alpha w_{50}/2 - w_{51}
  \big) u^8.
\end{align*}
It remains to consider necessity of the conditions of 
case e) for the lower bound $p =3$ of the interval $\Pi(7,3)$.
Let us assume that $\Lambda_7 P \neq 0$, but 
$\mu^{-3} P \in L^3(\Sigma,\mu)$.
Then, by \eqref{eq:ge} and \eqref{eq:T}, it holds
$\mu^{-3} \Lambda_7 P \in L^3(\Sigma,\mu)$,
implying that also $\mu^{-3} a_7$ and 
$\mu^{-3} d_7 \in L^3(\Sigma,\mu)$.
The function $d_7$ can be written in the form 
$d_7(u,v) = a u^3v^2 + b u^5 v + c u^7$ with certain 
coefficients $a,b,c$ so that, by Lemma~\ref{lem:g},
\[
 \infty > 
 \|\mu^{-3} d_7 \|_{3,\mu}^3 \ge  
 \int_\Delta \frac{|d_7|^3}{\mu^8} \ge 
 c_2 \int_\Delta 
 \frac{|a u^3 v^2 + b u^5 v + c u^7|^3}{(uv + u^3)^8}\, dudv.
\]
As in the proof of Lemma~\ref{lem:Lambda}, we use
the substitution $(u,v) = (st, s^2t(1-t))$ to obtain 
\[
 \int_\Delta 
 \frac{|a u^3 v^2 + b u^5 v + c u^7|^3}{(uv + u^3)^8}\, dudv \ge 
 \int_0^1 s^{-1}\, ds \cdot 
 \int_0^1 \big|a + (b-2a)t + (a-b+c)t^2\big|^3\, dt
 .
\]
Since the first integral on the right hand side is divergent,
we conclude that the second integral must vanish. This, in turn,
is possible only if $a=b=c=0$ and consequently $d_7=0$. An analogous 
argument shows that $a_7=0$. Hence, $\Lambda_7 P = 0$ in contradiction
to the assumption.
\end{proof}
Again, we state without proof that the conditions of case f) 
do not only guarantee 
membership in $W^{2,p}(\Omega)$ for all finite exponents $p$,
but even boundedness of the second partials,
$f \in W^{2,\infty}(\Omega)$.

In applications, the D-map $\bx$ can be assumed as given and 
fixed, while $z$ is to be determined according to the desired 
approximation purpose.
In a first step, the related D-map $\by = \Phi(\bx)$ in standard form 
is computed. Then, the relevant conditions on $w_{jk}$ can be 
transformed
to conditions on the sought coefficients $z_{jk}$ using the 
relations \eqref{eq:w}. These conditions 
are all linear so that the polynomials $z$ fulfilling 
them establish a linear subspace of functions 
$f = z \circ \bx^{-1} \in W^{2,p}(\Omega)$ for 
the targeted value of $p$. For instance, 
case c) yields
\[
 z_{10}=z_{01}=
 z_{11} = 0
 ,\quad 
 z_{21} = \alpha z_{20}+\beta z_{02}
 ,\quad 
 z_{12} = \gamma z_{20} + \delta z_{02}
 .
\]
According to Theorem~\ref{thm:C1}, these 
are just the conditions that guarantee $f$ to be 
continuously differentiable, and exactly these conditions 
are used for the experimental studies mentioned in 
the introduction. However, the theorem states 
that $f \not\in W^{2,2}(\Omega)$ unless the additional
conditions of case d) are satisfied. 
This indicates that the reported numerical results are baseless.

Further, our findings refute Theorem~1.4
in \cite{Reif:1994}, which states that the principal 
curvatures of D-patches are in $L^p$ for $p<4$. 
The crux is that the principal curvatures of the 
graph of a function $f$ are in $L^p$ if and only if 
$f$ is in $W^{2,p}$. Hence, in general, the 
principal curvatures of a D-patch are not in $L^p$ for $p \ge 2$.

\section{Conclusion}
\label{sec:conclusion}

We have derived a characterization of $W^{k,p}$-regularity for 
D-functions, as occurring in 
Isogeometric Analysis with degenerate geometry maps.
The results for $k=1$ show that $p=2$ comes for free
and that also larger values of $p$ can be obtained easily.
On the other hand, however, it is unclear why one should use at all 
a degenerate parametrization for the simulation of second order PDEs.
By contrast, the results for $k=2$ indicate that great care 
has to be taken in this case. In particular, the D-patch 
conditions known from the literature guarantee $C^1$, but not 
$W^{2,2}$. Further conditions must be incorporated 
to fix the approach. The following research tasks should be 
addressed in that respect: 

First, construct suitable bases of the subspaces fulfilling the 
$W^{2,p}$-conditions for $p \ge 2$. At least for $p=2$, 
in view of the relative simplicity of the conditions a) - d),
this should to be a solvable technical problem. Note that 
the more complicated conditions for large values of $p$
simplify considerably when dealing with low degree, say bicubic,
polynomials.

Second, scrutinize the approximation power of 
$W^{k,p}$-subspaces of D-functions with respect to 
different norms. Of course, numerical studies will be helpful,
but also an analytic treatment of the problem seems to be an appealing
challenge.

Third, derive similar regularity results for the important trivariate 
case. Given the fact that the present analysis of the bivariate 
case is not straightforward, severe technical difficulties must be 
expected with this project.


\bibliographystyle{alpha}
\bibliography{ref}

\end{document}